\documentclass[11pt]{article}
\usepackage{amssymb,amsmath,amscd}
\newtheorem{lem}{Lemma}[section]
\newtheorem{thm}{Theorem}[section]
\newtheorem{cor}{Corollary}[section]
\newtheorem{defn}{Definition}[section]
\newtheorem{conj}{Conjecture}

\newcommand{\f}[1]{\mathfrak{#1}}

\newcommand{\mb}{\mathbb}
\newcommand{\commentout}[1]{}

\newcommand{\sech}{{\rm sech}}
\newcommand{\mc}{\mathcal}
\newcommand{\arr}[1]{\left( \begin{array}{clcr} #1 \end{array} \right)}
\newcommand{\ol}[1]{\overline #1}
\newcommand{\xin}{Sp(n, \mb R)}
\newcommand{\meta}{Mp(n, \mb R)}
\newcommand{\diag}{{\rm diag}}
\begin{document}
\title{Functions on Symmetric Spaces and Oscillator Representation}
\author{Hongyu He \\
Department of Mathematics,\\
 Louisiana State University, \\
 Baton Rouge, LA 70803, U.S.A.\\
}
\date{}
\maketitle
\abstract{In this paper, we study the $L^2$ functions on $U(2n)/O(2n)$ and $\meta$. We relate them using the oscillator representation. We first study some isometries between various $L^2$ spaces using the compactification we defined in ~\cite{hhe99}. These isometries were first introduced by Betten-{\'O}lafsson in ~\cite{bo}~\footnote{I was informed by Prof. {\'O}lafsson of his work shortly after I finished this paper.}. We then give a description of the matrix coefficients of the oscillator representation $\omega$ in terms of algebraic functions on $U(2n)/O(2n)$. The structure of $L^2(U(2n)/O(2n))$ enables us to decompose the $L^2$ space of odd functions on $\meta$ into a finite {\it orthogonal direct sum}, from which  an orthogonal basis for $L^2(\meta)$ is obtained. In addition, our decomposition preserves both left and right $\meta$-action. Using this, we define the signature of tempered genuine representations of $\meta$. Our result implies that every genuine discrete series representation occurs as a subrepresentation in one and only one of
$(\otimes^p \omega) \otimes (\otimes^{2n+1-p} \omega^*)$ for $p$ with a fixed parity, generalizing some result in ~\cite{kv}. Consequently, we prove some results in the papers by Adam-Barbasch ~\cite{ab} and by Moeglin ~\cite{mo} without going through the details of the Langlands-Vogan parameter. In a weak sense, our paper also provides an analytic alternative to the Adam-Barbasch Theorem on Howe duality (~\cite{howe}).  }
\section{Introduction}
Harmonic analysis on symmetric spaces often involves analysis of spherical functions and spherical representations. In the compact case, Helgason's theorem gives the classification of spherical representations and the eigen-decomposition of the $L^2$ space. In the noncompact Riemanian case, the spherical functions are fairly complicated.
For $SL(2, \mb R)$ or more generally groups of real rank one, spherical functions are related to a class of special functions called the hypergeometric functions (~\cite{hel}, ~\cite{vi}). To study hypergeometric functions, typically one needs to study their series expansions and the differential equations defining them. Generally speaking, matrix coefficients of semisimple Lie groups, can be approached by  Harish-Chandra's power series expansion at \lq\lq infinity \rq\rq, or the Eisenstein integral (~\cite{hel} ~\cite{hel1}). Holomorphic methods also contribute to the studies of matrix coefficients.\\
\\
Oscillator representation first appeared in the papers of Shale, Segal and Weil as a \lq\lq projective \rq\rq representation of the symplectic group. It is quite different from the representations traditionally being studied by people like Gelfand and Harish-Chandra. Its construction is not as direct. In spite of this, the
 purpose of this paper is to show that,  for the oscillator representation and its tensor products, 
the matrix coefficients take a much simpler form and can be studied systematically. Furthermore, the matrix coefficients of the oscillator representation connects the analysis on $U(2n)/O(2n)$ with analysis on $\meta$. Our approach also allows us to analyze Howe duality in a way that is different from those of Moeglin, Adam-Barbasch (~\cite{howe}, ~\cite{mo}, ~\cite{ab}). \\
\\
Let $\mc F_n$ be the Fock space (~\cite{bar}). Let $\mc P_n$ be the space of polynomials in $n$ complex variables. Then $\mc P_n \subset \mc F_n$.
Let $(\omega, \mc F_n)$ be the Segal-Bargmann model of the oscillator representation of $\meta$ (~\cite{rr}, ~\cite{fo}). Put
\begin{equation}
\Lambda(g) = (\omega(g) 1(z), 1(z)),
\end{equation}
where $1(z)$ is the constant function $1$. 
Let $u, v \in \mc P_n$. Then $M_{\omega}(u \otimes v)(g)= (\omega(g)u, v)$ is called a matrix coefficient of $\omega$, with respect to $\mc P_n$. It can be easily shown that $\frac{M_{\omega}(u \otimes v)(g)}{\Lambda(g)}$ descends into a function of $\xin$. Throughout this paper, the function $M_{\omega}(u \otimes v)/\Lambda$ will always be regarded as a function on $\xin$. It turns out this function has intriguing algebraic and analytic properties.\\
\\
In ~\cite{hhe99}, we define an {\bf analytic compactifiction} $\mc H: \xin \rightarrow U(2n)/O(2n)$ using the Segal-Bargmann kernel. By ~\cite{he00}, this analytic compactification can be realized as follows.
Let $P_{2n}$ be the Siegel parabolic subgroup of $Sp(2n, \mb R)$. Let $X=Sp(2n, \mb R)/P_{2n}$. It is clear that $X \cong U(2n)/O(2n)$. Let $\xin \times \xin$ be diagonally embedded in $Sp(2n, \mb R)$. Then $\xin \times \xin$ acts on $X$. In ~\cite{he00}, we prove that there is a unique open and dense orbit $X_0$ in $X$ and furthermore $X_0 \cong \xin$. The identification of $X_0$ with the image of $\mc H$ is given in the Appendix of ~\cite{he00}. I should mention that our compactification
$\mc H$ is included in the list of a more general set of compactifications called {\it causal compactifications} (see Page 294 ~\cite{bo}). So the compacifications I defined in ~\cite{he00} have some overlap with the causal compactifications in Betten-{\'O}lafsson's list on Page 294 ~\cite{bo}.\\
\\
Let $f$ be a function on $\xin$. Let $f^0$ be the {\bf push-forward} of $f$ from $\xin$ to $U(2n)/O(2n)$.
In this paper, we prove two results concerning this push forward. 
\begin{thm}[Theorem 5.1 ~\cite{bo}]~\label{thm1}
The push forward $f \rightarrow f^0$ defines an isometry between $L^2(\xin, |\Lambda(g)|^{4n+2} d g)$ and $L^2(U(2n)/O(2n))$.
\end{thm}
 Here $|\Lambda(g)|$ is regarded as a function on $\xin$. This theorem is equivalent to Theorem 5.1 in ~\cite{bo} for the group $\xin$. The second result
states that the push-forward of $\frac{M_{\omega}(u \otimes v)(g)}{\Lambda(g)}$ is an algebraic function on $U(2n)/O(2n)$ and it is contained in $\oplus_{i \geq 0} C_{(0,0, \ldots,0, - 2i)}(U(2n)/O(2n))$ (see Theorem ~\ref{main}).\\
\\
Let us fix more notations before we state more results. Let $\epsilon$ be the nonidentity element in double covering of the identity in $\xin$. We call a representation $\pi$ of $\meta$, genuine if $\pi(\epsilon)=-1$. We say that a function on $\meta$ is odd if $f(\epsilon g)= - f(g)$.
Let $\omega_{p,q}=(\otimes^p \omega) \otimes (\otimes^q \omega^*)$. Let $\omega^c$ be the real representation $\omega$ equipped with a conjugate linear multiplication. Since $\omega$ is unitary, $\omega^* \cong \omega^c$.
Let $\mc P_{p,q}= (\otimes^p \mc P_n) \otimes (\otimes^q \mc P_n^c)$.  Using the compactification $\mc H$, we give a description of the matrix coefficients of $\omega_{p,q}$.
\begin{thm}~\label{thm2}
The push-forward of $(\omega_{p,q}(g) u, v)/\Lambda^p \overline{\Lambda}^q$ spans the space
$$\oplus_{\lambda \in 2 S_{q, p}} C_{\lambda}(U(2n)/O(2n)),$$
where $$S_{q, p}=\{(\lambda_1 \geq \lambda_2 \geq \lambda_q \geq 0 \ldots \geq 0 \geq \lambda_1^{\prime} \geq \ldots \geq \lambda^{\prime}_p) \mid \lambda_i, \lambda_j^{\prime} \in \mathbb Z \}$$
for $p+q \leq 2n$. Here $C_{\lambda}(U(2n)/O(2n))$ is a space of algebraic functions on U(2n)/O(2n) that is equivalent to the irreducible representation of $U(2n)$ with highest weight $\lambda$. See 3.4 for the definition of $S_{q,p}$ for $p+q \geq 2n+1$. 
\end{thm}
In connection with the eigen-decomposition of $L^2(U(2n)/O(2n))$, Theorem ~\ref{thm1} and Theorem ~\ref{thm2} imply the following.
\begin{thm}[Orthogonal Decompositions]
Let $\mc M_{p,q}$ be the linear span of matrix coefficients of $\omega_{p,q}$ with respect to $\mc P_{p,q}$.
 We have the following orthogonal decompositions
$$L^2_{-}(\meta)= \hat{\oplus}_{p \ odd} \mc M_{p, 2n+1-p}, \qquad L^2_{-}(\meta)= \hat{\oplus}_{p \ even} \mc M_{p, 2n+1-p}.$$
Furthermore, this decomposition preserves both the left and the right $\meta$ action.
\end{thm}
The space $\mc M_{0, 2n+1}$, a Hardy space, is studied in ~\cite{bo}.\\
\\
Let $\Pi^2_{-}(\meta)$ be the set of equivalence classes of irreducible tempered genuine representations of $\meta$.
Let $\Pi^2_{p,q}(\meta)$ be the set of equivalence classes of irreducible representations that are weakly contained  in $\omega_{p, 2n+1-p}$. In this paper, we prove that 
\begin{thm}~\label{thm3}
$$ \Pi^2_{p, 2n+1-p}(\meta) \cap \Pi^2_{p_1, 2n+1-p_1}(\meta) = \emptyset $$
if $p \equiv p_1 \pmod 2$ and $p \neq p_1$. Furthermore, 
$$\Pi^2_{-}(\meta)= \cup_{p \ odd} \Pi^2_{p, 2n+1-p}(\meta), \qquad \Pi^2_{-}(\meta)= \cup_{p \ even} \Pi^2_{p, 2n+1-p}(\meta).$$
\end{thm}
Let $\pi$ be a genuine tempered representation of $\meta$. Now we define the {\bf signature} of 
$\pi$ to be $(p^e, 2n+1-p^e, p^o, 2n+1-p^o)$ if 
$$ \pi \in \Pi^2_{p^e, 2n+1-p^e}(\meta) \cap \Pi^2_{p^o, 2n+1-p^o}(\meta),$$
where $p^e$ is even and $p^o$ is odd. Our definition of signature is closely related to the Howe duality ~\cite{howe}. It can be shown that a tempered genuine representation $\pi$ has signature $(p^e, 2n+1-p^e, p^o, 2n+1-p^o)$ if and only if $\pi$ occurs in $\mc R(\meta, \omega)$ for the dual pair $(\xin, O(p^e, 2n+1-p^e))$ and for the dual pair $(\xin, O(p^o, 2n+1- p^o))$. However, the connection is not obvious. Our $ \Pi^2_{p, 2n+1-p}(\meta)$ only includes irreducible representations that occur in the weak closure of $(\otimes^p \omega) \otimes (\otimes^{2n+1-p} (\omega^*))$, while $\mc R(\meta, \omega)$ refers to those irreducible representations that occur as quotients in $(\otimes^p \omega) \otimes (\otimes^{2n+1-p} (\omega^*))$ infinitesimally. \\
\\
In \cite{kv}, Kashiwara and Vergne proved that the representations occur in $\otimes^{2n+1} \omega$ are holomorphic discrete series representations of $\meta$. From Theorem ~\ref{thm3}, we obtain
\begin{cor} Fix a parity of $p$. Every irreducible genuine discrete series representation occurs as a subrepresentation of 
$$(\otimes^p \omega) \otimes (\otimes^{2n+1-p} (\omega^*))$$
for a unique $p$.
\end{cor}
The pair $(p, 2n+1-p)$ is exactly the signature.
For $p = 0, 2n+1$, the intertwining operator from $(\otimes^p \omega) \otimes (\otimes^{2n+1-p} \omega^*)$ to each discrete series representation can be computed and explored effectively. This has been done in ~\cite{kv}. For $p \neq 0, 2n+1$, the intertwining operator
$$(\otimes^p \omega) \otimes (\otimes^{2n+1-p} (\omega^*)) \rightarrow \pi$$
is not easy to describe. 
Our approach does not provide means to compute the signature of a discrete series representation. The question of obtaining the signature of a discrete series representation can be found in ~\cite{ab}. In fact, Adam and Barbasch proved 
\begin{thm}[Thm. 3.3, 5.1 ~\cite{ab}] Fix a parity of $p$. 
There is a one to one correspondence between the discrete series representations of $SO(p, 2n+1-p) (0 \leq p \leq 2n+1) $ and genuine discrete series representations of $\meta$. There is a one to one correspondence between the disjoint union of the admissible duals of $SO(p, 2n+1-p) (0 \leq p \leq 2n+1) $ and the genuine admissible dual  of $\meta$.
\end{thm}
In this paper, we prove the discrete series duality without going through the details of $K$-types. We also prove
\begin{thm}
Fix a parity of $p$. 
There is a one to one correspondence between the irreducible tempered representations of $SO(p, 2n+1-p) (0 \leq p \leq 2n+1) $ and irreducible genuine tempered  representations of $\meta$.
\end{thm}
In essence, our theorem is a $L^2$ Howe duality for $(O(p, 2n+1-p), \meta)$.

\section{Compactification of $Sp(n,\mb R)$}
The computations in this section are not new. They have been done by Betten-{\'O}lafsson  in a more general setting ~\cite{bo}.
\subsection{Setting}
Let $V$ be a $2n$-dimensional real vector space. Let $\{e_1, e_2, \ldots, e_{2n} \}$ be a 
basis for $V$. Let $J=\left( \begin{array}{clcr} 0 & I_n \\ -I_n & 0 \end{array} \right)$.
We equip $V$ with the following structures:
\begin{enumerate}
\item the standard symplectic form $\Omega(u,v)=u^t J v$;
\item the standard real inner product $(u,v)=u^t v$.
\item the complex structure
$$e_{n+i}= -i e_{i} \qquad (i=1,2, \ldots n).$$
\end{enumerate}
We denote $V$ equipped with the above complex structure by $V^{\mathbb C}$.
We equip $V^{\mathbb C}$ with the standard complex inner product $(,)^{\mathbb C}$.
Then $$(u,v)^{\mathbb C}=(u,v)+ i \Omega(u,v).$$
Let $Sp(n, \mb R)$ be the symplectic group preserving $\Omega(,)$. Let $O(2n)$ be the  orthogonal group preserving $(,)$.
Let $U(n)$ be the unitary group of $(V^{\mathbb C}, (,)^{\mathbb C})$.
Then $U(n)=Sp(n, \mb R) \cap
O(2n)$ and $U(n)$ is a maximal compact subgroup of $Sp(n,\mb R)$. \\
\\
Put
$$\f a= \{\diag(\lambda_1,\lambda_2,\ldots, \lambda_n,-\lambda_1,
-\lambda_2, \ldots, -\lambda_n) \mid \lambda_i \in \mb R, i \in [1,n]\}$$
Then $\f a$ is a maximal split Cartan subalgebra of $\f sp(n, \mb R)$. Let $A$ be the analytic group generated by
$\f a$. Let $K^0$  be the opposite group of $K$. Then $\xin$ has a $KAK$ decomposition and  $K \times K^o$ acts on $\xin$. \\
\\
Let $\mc S_{2n}$ be the space of $2n \times 2n$ symmetric unitary complex matrices. 
 $\mc S_{2n}$ can be identified with $U(2n)/O(2n)$ as follows
 \begin{equation}~\label{mapj}
 j: g \in U(2n) \rightarrow g g^t \in \mc S_{2n}.
 \end{equation} 
 The group $U(n) \times U(n)$ acts on $\mc S_{2n}$ as follows:
$$\tau(k_1,k_2) s = \diag(k_1,k_2) \  s  \ \diag(k_1^t, k_2^t) \qquad (k_1,k_2 \in U(n)).$$
\begin{thm}[Analytic Compactification of $\xin$, ~\cite{hhe99}]~\label{com} Let $K=U(n)$. For $g=k_1 \exp H k_2 \in \xin$,  let
$$\mc H(k_1 \exp H k_2)=\left( \begin{array}{clcr} \ol{k_1} & 0 \\
0 & \ol{k_2}^t \end{array} \right) 
\left( \begin{array}{clcr} \tanh H^{\mb C} & -i \sech H^{\mb C} \\
-i \sech H^{\mb C} & \tanh H^{\mb C} \end{array} \right)
\left( \begin{array}{clcr} k_1^{-1} & 0 \\ 0 & \ol{k_2} \end{array} \right).
$$
where $H=(\lambda_1,\lambda_2,\ldots, \lambda_n,-\lambda_1,-\lambda_2,\ldots,-\lambda_n)$ and 
$H^{\mb C}=(\lambda_1,\lambda_2,\ldots, \lambda_n)$.
Then
\begin{enumerate}
\item $\mathcal H$ is a well-defined map from $\xin$ to $\mathcal S_{2n}$;
\item $\mc H$ is an analytic embedding;
\item the image of $\mc H$ is open and dense in $\mc S_{2n}$;
\item  Identify
$K \times K^o$ with $K \times K$ by
$$(k_1, k_2) \rightarrow (\ol{k_1}, k_2^{-1}).$$
Then $\mc H: \xin \rightarrow \mathcal {S}_{2n}$ is $K \times K^o$-equivariant;
\item Let $f$ be a matrix coefficient of a unitary irreducible representation
of $\xin$. Then $f$ can be pushed forward to a continuous
function $f^0$ on $\mc S_{2n}$.
\end{enumerate}
\end{thm}
In short, $\mc H$ is an analytic compactification
of $\xin$.
 For an arbitrary function $f$ on $\xin$, we will use
$f^0$ to denote the push-forward of $f$ to $\mc H(\xin)$. If
$f^0$ has a continuous extension over $\mc S_{2n}$, then the extension must be unique.
In such a case, $f^0$ will
be used to denote the unique extension. See ~\cite{hhe99} for details.

\subsection{Invariant Measure on $U(2n)$ with respect to the Generalized Cartan Decomposition}
Generalized Cartan Decomposition is studied in Heckman and Schlichtkrull's book ~\cite{hs}.
Unlike the Cartan decomposition which is associated with a symmetric pair $(G, G^{\sigma})$,
generalized Cartan decomposition is built on a symmetric triple $(G, G^{\sigma}, G^{\tau})$ with
$(\sigma,\tau)$ a commuting pair of involutions. We start with the structure theory. \\
\\
Define
$$\sigma(g)=\arr{I_n & 0_n \\ 0_n & -I_n  } g \arr{I_n & 0_n \\ 0_n & -I_n }; \tau(g)=\overline{g}.$$
Let  $U=U(2n), K=O(2n)$ and $H=U(n) \times U(n)$. 
Let
$$\f p= \{ i B \mid B^t=B, B \in \f {gl}(2n, \mb R) \}$$
and
$$\f q= \{ \arr{ 0 & A \\ -\overline{A}^t & 0 } \mid A \in \f {gl}(n,\mathbb C) \}.$$
Then 
$$U^{\tau}=K, U^{\sigma}=H, \f k \oplus \f p= \f u, \mathfrak h \oplus \mathfrak q =\mathfrak u.$$
Let $T_{\f pq}$ be the torus
consisting of
$$t=\arr{\cos \theta & -i \sin \theta \\ -i \sin \theta & \cos \theta}$$
where $\theta \in \mathbb R^n$. $\f t_{\f pq}$ is a maximal Abelian subalgebra in $\f p \cap \f q$. Then 
 the generalized Cartan
decomposition (Theorem 2.6, P. 194, ~\cite{hs}) says that every $g \in U$ can be written as
a product $k(g) t(g) h(g)$. Moreover,
$$d g= J(t) d t d k d h$$
with
$$J(t)= \prod_{\alpha \in \Sigma_{\f pq}^+} | t^{\alpha}- t^{-\alpha}|^{m_{\alpha}^+}
|t^{\alpha}+ t^{-\alpha}|^{m_{\alpha}^-}.$$
Here $m_{\alpha}^+= \dim {\f u}_{\alpha}^+$ and $m_{\alpha}^-=\dim {\f u}_{\alpha}^-$
are the multiplicities of the complex root space
$$\f u_{\alpha}^+=\{ x \in \f {u}_{\mathbb C} \mid \sigma \tau(x)=x, [h, x]= \alpha(h) x , \forall h \in \f t_{\f pq}  
 \} $$
$$\f u_{\alpha}^-=\{ x \in \f {u}_{\mathbb C} \mid \sigma \tau(x)=-x, [h, x]= \alpha(h) x , \forall h \in \f t_{\f pq}  
 \}. $$
This result is due to Heckman-Schlichtkrull.
\begin{thm} For $$t=\arr{\cos \theta & -i \sin \theta \\ -i \sin \theta & \cos \theta},$$ 
$$J(t)= c | \prod_{i > j} (\cos 4 \theta_i - \cos 4 \theta_j)| | 
\prod_{i} \cos 2 \theta_i |.$$
Here $c$ is some constant.
\end{thm}
Proof: Notice
$$\sigma \tau \arr{U_1 & V \\ - \overline{V}^t & U_2 }=\arr{\overline{U_1} & - \overline{V} \\ V^t & \overline{U_2} }$$
Thus
$$\f u^+=\{\arr{{U_1} & i V \\ i V^t & {U_2} } \mid  U_1^t+U_1=0=U_2^t+U_2, \ \ 
U_1, U_2, V \in \f {gl}(n, \mb R) \}$$
$$\f u^+_{\mathbb C} =\{\arr{{U_1} &   V \\  V^t & {U_2} } \mid  U_1^t+U_1=0=U_2^t+U_2, \ \ 
U_1, U_2, V \in \f {gl}(n, \mb C) \}$$
Similarly,
$$\f u^-=\{\arr{{i U_1} &  { V} \\  -V^t & {i U_2} } \mid  U_1^t=U_1, U_2^t=U_2, \ \ 
U_1, U_2, V \in \f {gl}(n, \mb R) \}$$
$$\f u^-_{\mathbb C}=\{\arr{{ U_1} & { V} \\  -V^t & { U_2} } \mid  U_1^t=U_1, U_2^t=U_2, \ \ 
U_1, U_2, V \in \f {gl}(n, \mb C) \}.$$
Let $C=\frac{1}{\sqrt{2}} \arr{I_n & I_n \\ I_n & -I_n }$. Then
$$C \arr{ 0 & i \theta \\ i \theta & 0} C^{-1}= i \arr{\theta & 0 \\ 0 & -\theta}$$
$$C \arr{U_1 & V \\ W & U_2 } C^{-1}=\frac{1}{2}\arr{U_1+U_2 + W+V & U_1-U_2+W-V \\
U_1-U_2-W+V & U_1+U_2 -W-V }$$
Therefore 
$$C \f u^+_{\mathbb C} C^{-1}=\{ \arr{ X & Y \\ Z & -X^t } \mid Y^t+Y=Z^t+Z=0 \}$$
$$C \f u^-_{\mathbb C} C^{-1}=\{ \arr{ X & Y \\ Z & -X^t } \mid Y^t=Y,Z^t=Z \}$$
Let $e_i(\theta)=\theta_i$. Then 
$$m_{\pm e_i \pm e_j}^+=m_{\pm e_i \pm e_j}^-=1 \qquad (i \neq j) $$
$$m_{2 e_i}^+=0; m_{2 e_j}^-=1.$$
By Theorem 2.7, P. 194, ~\cite{hs}, 
\begin{equation}
\begin{split}
J(t)=  & | 2^{2n^2-n} \prod_{i > j} 
\cos(\theta_i-\theta_j)\sin(\theta_i-\theta_j) \cos(\theta_i+\theta_j) \sin(\theta_i+\theta_j) \prod_i 
\cos(2 \theta_i) | \\
=  &  2^{n^2} | \prod_{i > j} \sin(2 \theta_i -2 \theta_j) \sin (2 \theta_i+ 2 \theta_j) 
\prod_{i} \cos( 2\theta_i) | \\
= & c |\prod_{i > j} \cos(4 \theta_i)-\cos( 4 \theta_j) \prod_i \cos (2 \theta_i)|
\end{split}
\end{equation}
Q.E.D. \\
\\
Under the generalized Cartan decomposition, the map $j: U(2n) \rightarrow U(2n)/O(2n)$ (see Equ. ~\ref{mapj} ) becomes,
$$ h \arr{\cos \theta & -i \sin \theta \\ -i \sin \theta & \cos \theta} k
\rightarrow h \arr{\cos 2 \theta & -i \sin 2 \theta \\ -i \sin 2 \theta & \cos 2 \theta}  
h^t \in \mc S_{2n} \qquad (h \in H, k \in K).$$
Thus the invariant measure $ds$ for $\mc S_{2n}$  is given by
$J(t) d h d \theta$ for
$s=h t^2 h^t$.
\begin{cor}~\label{measureons} Under the generalized Cartan decomposition 
$$s=h \arr{\cos  \theta & -i \sin  \theta \\ -i \sin  \theta & \cos \theta}  h^t,$$
the $U(2n)$ invariant measure on $\mc S_{2n}$ is given by
$$d s=c |\prod_{i > j} (\cos(2 \theta_i)-\cos( 2 \theta_j)) | \ | \prod_i \cos  \theta_i | d h d \theta.$$
\end{cor}

\subsection{The Jacobian $| \frac{ d \mc H(g)} { d g}|$ }
Consider the group $G=\xin$ and
the root system
$$\Sigma(\f a, \f g)=\{ \pm e_i \pm e_j (i > j); \pm 2 e_i (i, j \in [1, n]) \}.$$
Under the Cartan decomposition $g= k_1 \exp H k_2$ with $(k_1, k_2 ) \in U(n) \times U(n)$,
\begin{equation}~\label{measureonxin}
\begin{split}
  d g = & |\prod_{\alpha \in \Sigma^+(\f a, \f g)} \sinh \alpha(H) | 
  d k_1 d k_2  d \lambda \\ 
=& | \prod_{i > j}  \sinh (\lambda_i-\lambda_j) \sinh (\lambda_i+ \lambda_j) 
\prod_i \sinh 2\lambda_i| d \lambda d k_1 d k_2 \\
=& 2^{-n^2} \prod_{i > j} | (\exp (\lambda_i-\lambda_j)-\exp(-\lambda_i+\lambda_j))
(\exp (\lambda_i+ \lambda_j)-\exp(-\lambda_i-\lambda_j))| \\
& \prod_i |(\exp (2 \lambda_i)-
\exp (- 2\lambda_i))|  d k_1 d k_2 d \lambda \\ 
= &  2^{-n^2} \prod_{i > j}  |(\exp 2 \lambda_i+ \exp (-2 \lambda_i)- 
\exp 2 \lambda_j - \exp (-2 \lambda_j))| \\
&  \prod_i |(\exp ( \lambda_i)-
\exp (- \lambda_i))(\exp(\lambda_i)+ \exp(-\lambda_i))|  d k_1 d k_2  d \lambda
\end{split}
\end{equation}
Recall from Theorem ~\ref{com}, $\mc H$ is $U(n) \times U(n)$ equivariant and
$$\mc H: \exp H=\diag( \exp \lambda, \exp(-\lambda)) 
\rightarrow \arr{\tanh \lambda & - i \sech \lambda \\
- i \sech \lambda & \tanh \lambda }
=\arr{\cos  \theta & -i \sin  \theta \\ -i \sin  \theta & \cos \theta}. $$
Thus
$$\cos \theta= \tanh \lambda; \qquad \sin \theta= \sech \lambda. $$
It follows that
$$- \sin \theta_i d \theta_i = \sech^2 \lambda_i d \lambda_i; \qquad \cos \theta_i d \theta_i= -\sech \lambda_i \tanh \lambda_i d \lambda_i.$$
Consequently,
$$\frac{d \theta_i}{ d \lambda_i}=- \sech \lambda_i = - \sin \theta_i \qquad
\frac{ d \lambda_i} { d \theta_i}= -\csc \theta_i= -  \cosh \lambda_i .$$
We state the following theorem concerning the Jacobian of $\mc H$.
\begin{thm}
Let $d s$ be a $U(2n)$-invariant measure on $\mc S_{2n}$. Let $ d g$ be a Haar measure
on $\xin$. Let $g= k_1 \exp H k_2$ with $H=\diag(\lambda_1, \ldots, \lambda_n, 
- \lambda_1, \ldots,-\lambda_n)$. Then
$${ d \mc H(g)}=c (\prod_{i} \sech \lambda_i)^{2n+1} d g.$$
\end{thm}
Proof: For every $(i > j)$, 
\begin{equation}
\begin{split}
& (\exp(2 \lambda_i)+ \exp(-2 \lambda_i)- \exp (2 \lambda_j)- \exp (-2 \lambda_j))
(\sech^2 \lambda_i \sech^2 \lambda_j) \\
 = & 16 \frac{\exp(2 \lambda_i)+ \exp(-2 \lambda_i)- \exp (2 \lambda_j)- \exp (-2 \lambda_j)}
{(\exp \lambda_i+ \exp (-\lambda_i))^2(\exp \lambda_j+ \exp (- \lambda_j))^2} \\
= & 16 \frac{(\exp \lambda_i+ \exp (-\lambda_i))^2 - (\exp \lambda_j + \exp (-\lambda_j))^2}
{(\exp \lambda_i+ \exp (-\lambda_i))^2(\exp \lambda_j+ \exp (- \lambda_j))^2} \\
= & 16 (\frac{1}{(\exp \lambda_j -\exp (-\lambda_j))^2}-\frac{1}{(\exp \lambda_i -\exp (-\lambda_i))^2}) \\
= & 4 (\sech^2 \lambda_j - \sech^2 \lambda_i) \\
= & 4 (\sin^2 \theta_j -\sin^2 \theta_i) \\
= & 2( \cos 2 \theta_i -\cos 2 \theta_j)
\end{split}
\end{equation}
For every $i$, we have
\begin{equation}
 \sech^2 \lambda_i ( \exp \lambda_i- \exp (- \lambda_i)) (\exp \lambda_i + \exp (-\lambda_i))
= 4 \tanh \lambda_i = \cos \theta_i 
\end{equation}
\begin{equation}
-\sech \lambda_i d \lambda_i = d \theta_i 
\end{equation}
Multiplying these three equations together and taking the absolute value, we obtain
\begin{equation}
\begin{split}
 & \prod_{i > j} |(\exp(2 \lambda_i)+ \exp(-2 \lambda_i)- \exp (2 \lambda_j)- \exp (-2 \lambda_j))
(\sech^2 \lambda_i \sech^2 \lambda_j) | \\ 
& \prod_i |\sech^2 \lambda_i ( \exp \lambda_i- \exp (- \lambda_i)) (\exp \lambda_i + \exp (-\lambda_i)) |
\prod_i |\sech \lambda_i|  d \lambda_i \\
= & c \prod_{i > j} |( \cos 2 \theta_i -\cos 2 \theta_j) | \prod_i |\cos \theta_i | d \theta_i 
\end{split}
\end{equation}
Since $\mc H$ is $U(n) \times U(n)$-equivariant,
from Cor ~\ref{measureons} and Equation ~\ref{measureonxin}, we obtain
$$ c (\prod_i \sech \lambda_i)^{2n+1} d g = d \mc H(g).$$
Here $c$ is used as a symbolic constant.
Q.E.D.
\begin{cor}[See Theorem 5.1 ~\cite{bo}]
The push forward $f \rightarrow f^0$ defines an isometry between $L^2(\xin, (\prod \sech \lambda_i(g))^{2n+1} d g)$ and $L^2(\mc S_{2n}, d s)$.
\end{cor}

\section{Functions on $\mc S_{2n}$ and $\xin$}

\subsection{Helgason's Theorems}
Let $(U,K)$ be a reductive symmetric pair of compact type (see ~\cite{hel}).
Let $\pi$ be
an irreducible representation of $U$. $\pi$ is said to be spherical if there exists
a nonzero vector that is fixed by
$\pi(K)$.
 Let $\f p$ be the orthogonal complement of $\f k$ in 
$\f u$. Let $\f t_{\f p}$ be a maximal abelian subalgebra of $\f p$.
Let $M$ be the centralizer of $\f t_{\f p}$ in $K$. Let $\f t$ be a maximal toral
subalgebra of $\f u$ containing $\f t_{\f p}$.
\begin{thm}[Helgason's Theorem]
Let $\pi$ be an irreducible representation of $U$. Then $\pi$ is spherical if and only
if $\pi(M)$ leaves the highest weight vector fixed. Furthermore the spherical vector
is unique up to a scalar.
\end{thm}
 Let $\mathcal S=U/K$. 
 An immediate consequence of
 Helgason's theorem is the
decomposition theorem of square integrable functions on $\mc S$. 
\begin{thm} Let $L^2(\mc S)$ be the space of 
square integrable
functions with respect to the $U$-invariant measure on $\mc S$. Then
 $L^2(\mc S)$ is the closure of
$$\oplus_{\mu \ spherical} C_{\mu}(\mc S)$$
Here $C_{\mu}(\mc S)$ is equivalent to the irreducible unitary
representation with highest weight $\mu$.
\end{thm}
Now let $U=U(n)$, $K=O(n)$. Then $(U,K)$ is a reductive symmetric pair of compact type. Even though Helgason's original theorem assumes that $U$ is simply connected and semisimple,
 it remains valid  for $(U(n), O(n))$. Let
$$\f t_{\f p}=\{\diag(i \theta_1,i\theta_2, \ldots i\theta_{n}) \mid \theta_i \in \mathbb R\}.$$
Then the Weyl group $W(U,K)$ is the permutation group on $\{\theta_i \}$. The centralizer
$$M=\{\epsilon=diag(\epsilon_1, \epsilon_2, \ldots, \epsilon_{n}) \mid \epsilon_i = \pm 1\} 
\subseteq T_{\f p}.$$
Let $V_{\lambda}$ be an irreducible representation of
$U(n)$ with the highest weight 
$$\lambda=(\lambda_1 \geq \lambda_2 \geq \ldots \geq \lambda_{n}).$$
Notice that for $v_0$, the highest weight vector of $(\pi,V_{\lambda})$,
$$\pi(\epsilon) v_0= \Pi_{i=1}^{n} (\epsilon_i)^{\lambda_i} v_0.$$
Therefore
$V_{\lambda}$ has a $K$-fixed vector if and only if every $\lambda_i$ is {\em even}. We say that
$\lambda$ is {\bf even} if every $\lambda_i$ is even.
We obtain
\begin{cor} ~\label{un} Consider the reductive symmetric pair $(U(n), O(n))$. An irreducible representation $V_{\lambda}$ of $U(n)$ is spherical
if and only
if $\lambda$ is {\em even}.
Moreover, we have
$$L^2({\mathcal S})= \hat{\oplus}_{\lambda \ \ even} C_{\lambda}(\mathcal S)$$
where $C_{\lambda}(\mathcal S)$ is an irreducible representation with highest weight 
$$\lambda=(\lambda_1 \geq \lambda_2 \ldots \geq \lambda_{n}).$$
\end{cor}

\subsection{The Metaplectic Function}
For a $2n \times 2n$ matrix $g$, define $$C_g=\frac{1}{2}(g -JgJ).$$ Suppose
that
$$g = \arr{ A & B \\ C & D }.$$ Then
$$C_g=\arr{\frac{A+D}{2} & \frac{B-C}{2} \\ \frac{C-B}{2} & \frac{A+D}{2}}.$$
Define
$$C_g^{\mb C}=\frac{A+D}{2}+i\frac{B-C}{2}.$$
It is known that $C_g^{\mb C} \in GL_n(\mb C)$ when $g \in \xin$ (See ~\cite{rr}).
One can now write down the metaplectic group $\meta$ precisely as follows (see ~\cite{rr})
$$\meta=\{(\lambda,g) \mid g \in \xin, \lambda^2 \det (C_g^{\mb C})=1\}.$$
Define the {\bf metaplectic function} $\Lambda$ on $\meta$ to be
$\Lambda(\lambda, g)=\lambda $. Then $\Lambda^2=\det (C_g^{\mathbb C})^{-1}$ is a function on 
$\xin$. $|\Lambda(\lambda, g)= |\det C_g^{\mb C}|^{-\frac{1}{2}}$ is also a function on $\xin$. \\
\\
Let $g=k_1 \exp H k_2$ be a $KAK$ decomposition of $\xin$. Recall from ~\cite{hhe99} that
$$C_g^{\mb C}=k_1 \cosh (H^{\mb C}) k_2$$
Thus 
$$\Lambda^2(g)= \det (C_g^{\mb C})^{-1}= \det (k_2^{-1} \sech (H^{\mb C}) k_1^{-1})=\det (\overline{k_2} \sech (H^{\mb C}) 
\overline{k_1})$$
and
$|\Lambda^2(g)|= \prod_{i} \sech \lambda_i(g)$.
Over $\mathcal S_{2n}$,
define a function
$${\det}_{12}:  s=\left( \begin{array}{clcr} s_{11} &  s_{12}
\\  s_{12}^t & s_{22} \end{array} \right) \in \mc S_{2n} \rightarrow \det s_{12}.$$
From our definition of $\mc H$, we see immediately that
$$\Lambda(g)^2=\det (C_g^{\mb C})^{-1} =\det (\overline{k_2} \sech (H^{\mb C}) 
\overline{k_1})=  {\det}_{12}( i \mc H(g)).$$
So the push-forward of $\Lambda^2$ under $\mc H$ is
$$(\Lambda^2)^0=i^n {\det}_{12}. $$
Let $V$ be the standard $2n$-dimensional Hilbert space with orthonormal basis
$$\{ e_1, e_2, \ldots e_{2n} \}. $$
Let $U(2n)$ act on $V$ canonically as unitary operators.
Then $(\wedge^n, \wedge^n V)$ becomes a unitary representation of $U(2n)$.
In particular,
$${\det}_{12}(s)=(\wedge^n(s) e_1 \wedge e_2 \wedge \ldots \wedge e_n, e_{n+1} \wedge e_{n+1}
\wedge \ldots \wedge e_{2n}).$$
 We summarize
our discussion in the following theorem.
\begin{thm}~\label{phi12} The push-forward of $\Lambda^2$ under $\mc H$ is
$$i^n {\det}_{12}(s)= i^n (\wedge^n(s)  e_1 \wedge e_2 \wedge \ldots \wedge e_n, e_{n+1} \wedge e_{n+1}
\wedge \ldots \wedge e_{2n}).$$
where  $e_1 \wedge e_2 \wedge \ldots \wedge e_n$ is a highest weight vector of $\wedge^n$ and
$e_{n+1} \wedge e_{n+2} \wedge \ldots \wedge e_{2n}$ is a lowest weight vector.
\end{thm}

\subsection{An Isometry Between $L^2(\mc S_{2n})$ and $L^2_{-}(\meta)$}
Let $\epsilon$ be the nonidentity element in the metaplectic lifting of the identity
element in $\xin$. We say that a function $f$ on $\meta$ is odd  if
$$ f(\epsilon \tilde g) = - f(\tilde g) \qquad (\tilde g \in \meta);$$
we say that a function $f$ on $Mp_{2n}(\mb R)$ is even  if
$$ f(\epsilon \tilde g) =  f(\tilde g) \qquad (\tilde g \in \meta).$$
Let $L^2_{-}(\meta)$ be the space of odd square integrable functions of $\meta$.
Let $L^2_{+}(\meta)$ be the space of even square integrable functions of $\meta$.
We identify $L^2_{+}(\meta)$ with $L^2(\xin$.
For any $ f \in L^2_{-}(\meta)$, observe that
$f \Lambda^{2n+1}$ is an even function. Regarding $ f \Lambda^{2n+1}$ as a function on $\xin$, 
we define
$$I: 
  f \in L^2_{-}(\meta) \rightarrow (f \Lambda^{-2n-1})^0.$$
$I(f)$ is a function on $\mc S_{2n}$.
\begin{thm}[See Theorem 7.1 ~\cite{bo}] ~\label{iso}
With a proper choice of the invariant measure on $\meta$, the map $I$ defines isometry from $L^2_{-}(\meta)$ to $L^2(\mc S_{2n})$.
\end{thm}
Proof: Let $f, h \in L^2_{-}(\meta)$.
Then
\begin{equation}
\begin{split}
(f, h)_{L^2(\meta)}= & \int_{\tilde g \in \meta} f(\tilde g) \overline{h(\tilde g)} d \tilde g \\
 = &  \int_{\meta} (f \Lambda^{-2n-1})(\tilde g) (\overline{h \Lambda^{-2n-1}})(\tilde g)
 |\Lambda^{4n+2}(\tilde g) d \tilde g \\
 = & 2 \int_{\xin} (f \Lambda^{-2n-1})(g) (\overline{h \Lambda^{-2n-1}})( g) 
  \prod_i \sech^{2n+1} \lambda_i(g) d g \\
 = & 2 \int_{\mc S_{2n}} I(f)(s) \overline{I(h)}(s) d s \\
 =  & 2 (I(f), I(h))_{L^2(\mc S_{2n})}.
 \end{split}
 \end{equation}
 Q.E.D.

\subsection{An Orthogonal Decomposition of $L^2(\mathcal S_{2n})$}
Let $\mathbb Z^n$ be the integral lattice of $n$ dimension. Let $2 \mathbb Z^{n}$ be the sublattice of even integers.
Let $\det$ be the determinant function of $s \in \mathcal S_{2n}$. The determinant function is a weight $(-2,-2, \ldots, -2)$ function on $\mathcal S_{2n} \cong U(2n)/O(2n)$. For $p+q \leq 2n$, we define a subset of $\mathbb Z^{2n}$:
$$S_{p,q}= \{ \overbrace{\lambda_1 \geq \lambda_2 \geq \ldots \geq \lambda_p}^p \geq \overbrace{ 0 =0 \ldots =0}^{2n-p-q} \geq \overbrace{\mu_1 \geq \ldots \geq \mu_q}^q \}.$$
For $p+q \geq 2n+1$, we define to $S_{p,q}$ be the union of $S_{s,t}$ with $s \leq p, t \leq q$ and $s+t \leq 2n$.
Write
$$\mathbf m=(\overbrace{m,m, \ldots m}^{2n});$$
$$S_{p,q}+ \mathbf m= \{ \lambda+ \mathbf m \mid \lambda \in S_{p,q} \}.$$
Let $\mc O(\mc S_{2n})$ be the space of regular functions on $\mc S_{2n}$.
\begin{thm}~\label{deco}
$$\mathcal O(\mathcal S_{2n})= \oplus_{p \ even, p+q=2n+1} C_{ \lambda \in 
S_{p,q}+ \mathbf{p} \cap 2 \mathbb Z^{2n} }(\mathcal S_{2n}).$$
$$\mathcal O(\mathcal S_{2n})= \oplus_{p \ even, p+q=2n+1} C_{ \lambda \in 
S_{p,q}- \mathbf{p} \cap 2 \mathbb Z^{2n} }(\mathcal S_{2n}).$$
$$\mathcal O(\mathcal S_{2n})= \oplus_{q \ even, p+q=2n+1} C_{ \lambda \in 
S_{p,q}- \mathbf{q} \cap 2 \mathbb Z^{2n}}(\mathcal S_{2n}).$$
$$\mathcal O(\mathcal S_{2n})= \oplus_{q \ even, p+q=2n+1} C_{ \lambda \in 
S_{p,q} + \mathbf{q} \cap 2 \mathbb Z^{2n}}(\mathcal S_{2n}).$$
\end{thm}
Proof: We will only prove the first statement. It suffices to show that for any $\lambda$ even, there exists a unique $p+q=2n+1$ with $p$ even such that $\lambda \in S_{p,q}+ \mathbf{\frac{p}}$. Observe that  $\lambda \in S_{2i, 2n+1-2i}+\mathbf 2i$ if and only if
$\lambda_{2i+1} \leq 2i$ and $\lambda_{2i-1} \geq 2i$. 
If $\lambda_{1} \leq 0$, then $\lambda \in S_{0,2n+1}$. Otherwise, $\lambda_{1} \geq 2$.  If $\lambda_3 \leq 2$, then
$\lambda \in S_{2, 2n-1}+ \mathbf 2$. Otherwise $\lambda_3 \geq 4$. If $\lambda_5 \leq 4$, then $\lambda \in S_{4, 2n-3}+\mathbf 4$. We continue on this process. If $\lambda_{2i+1} \leq 2i$ and $\lambda_{2i-1} \geq 2i$, then $\lambda \in S_{2i, 2n+1-2i}+\mathbf 2i$. Finally, if $\lambda_{2n-1} \geq 2n$, $\lambda \in S_{2n,1}+\mathbf 2n$. Q.E.D.
\begin{thm}
$$\mathcal O(\mathcal S_{2n})= \oplus_{p+q=2n} C_{\lambda \ even, \lambda \in 
S_{p,q}+ \mathbf{2p}}(\mathcal S_{2n}).$$
\end{thm}
Proof: Suppose that $\lambda$ is even. Notice that $\lambda \in 
S_{p,q}+ \mathbf{2p}$ if and only if $\lambda_p \geq 2p \geq \lambda_{p+1}$. If $\lambda_1 \leq 0$, then
$\lambda \in S_{0, 2n}$. Otherwise, $\lambda_1 \geq 2$. If $\lambda_2 \leq 2$, then $\lambda_p \geq 2p \geq \lambda_{p+1}$ for $p=1$. So $\lambda \in S_{1, 2n-1}+ \mathbf 2$. Otherwise, $\lambda_2 \geq 4$. If $\lambda_3 \leq 4$, then  $\lambda_p \geq 2p \geq \lambda_{p+1}$ for $p=2$. So $\lambda \in S_{2, 2n-2}+ \mathbf 4$. We continue on this process until $p=2n$ for which $\lambda_{2n} \geq 4n$. So $\lambda \in S_{2n,0}+\mathbf 4n$. Q.E.D.\\
\\
Let $F$ be an $L^2$ function on $\mathcal S_{2n}$. Then $F$ has a Helgason-Peter-Weyl expansion
$$F = \sum_{\lambda even} F_{\lambda}.$$
where $F_{\lambda} \in C_{\lambda}(\mathcal S_{2n})$.
 Let $S$ be a subset of $\mathbb Z^{2n}$. Let
$L^2(\mathcal S_{2n})(S)$ be the subspace of $L^2$-functions whose HPW expansion only contains
$F_{\lambda}$ with $\lambda \in S$. We have
\begin{thm}
$$L^2(\mathcal S_{2n})=\oplus_{p \ even, p+q=2n+1}  L^2(\mathcal S_{2n})( 
S_{p,q}+ \mathbf{p} \cap 2 \mathbb Z^{2n});$$
$$L^2(\mathcal S_{2n})= \oplus_{q \ even, p+q=2n+1} L^2(\mathcal S_{2n})( 
S_{p,q}- \mathbf{q} \cap 2 \mathbb Z^{2n});$$
$$L^2(\mathcal S_{2n})=  \oplus_{p+q=2n} L^2(\mathcal S_{2n})(S_{p,q}+ 2 \mathbf p \cap 2 \mathbb Z^{2n}).$$
\end{thm}

\section{Matrix Coefficients of The Oscillator Representation}
Let me first introduce the Bargmann-Segal Model of the oscillator representation.
Let $V=\mb C^n$ and $$d \mu(z)= \exp (-\frac{1}{2} \|z\|^2) d z_1 \ldots dz_n
$$ be the Gaussian measure. For simplicity, 
let $ d z$ be the Euclidean measure on $V$. Let $\mc P_n$ be the
space of polynomials on $V$. We define an inner product on $\mc P_n$
$$(f,g)=\int f \ol g d \mu(z).$$
It is well-known that
$$(z^{\alpha}, z^{\beta})=0 \qquad (\alpha \neq \beta)$$
$$(z^{\alpha}, z^{\alpha})=2^{\alpha} \alpha !$$
(See ~\cite{bar}). Here we follow the multi-index convention and
$$2^{\alpha}=2^{|\alpha|}, \qquad |\alpha|=\sum_{i=1}^n \alpha_i.$$
Now let $\mc F_n$ be the completion of $\mc P_n$. Then $\mc F_n$ is
precisely the space of square integrable analytic functions on $V$
(see ~\cite{bar}). \\
\\
For $(\lambda, g) \in \meta$, we define
$$\omega(\lambda,g) f(z) =\int_V \lambda \exp (\frac{1}{4} (iz^t, \ol w^t)
\mc H(g) \arr{iz \\ \ol w}) f(w) d \mu(w)$$
Then $(\omega, \mc F_n)$ is a faithful unitary representation of $\meta$ (see ~\cite{rr}).
This model is often called the {\it Bargmann-Segal} model of the
oscillator representation. 
\subsection{Matrix coefficient Maps}

Let $\phi,\psi \in \mc P_n$. Define {\bf the map of matrix coefficient}
$$M_{\omega}: \mc P_n \otimes \mc P_n \rightarrow C^{\infty}(\meta)$$
as follows
$$M_{\omega}(\phi \otimes \psi)(g)=(\omega(g) \phi, \psi), \qquad \forall \ g \in \meta.$$
We have
\begin{equation}
\begin{split}
M_{\omega}({\phi \otimes \psi})(\lambda, g) =&(\omega(\lambda,g) \phi(z), \psi(z)) \\
=&\int_{V \times V} \lambda \exp (\frac{1}{4} (i z^t, \ol w^t) \mc H(g) 
\arr{iz \\ \ol w})
 \phi(w) \overline{\psi(z)} d \mu(w) d \mu(z) \\
=& \int_{V \times V} \lambda \exp (\frac{1}{4} (\ol z^t,  \ol w^t) \mc H(g) 
\arr{\ol z \\  \ol w})
\phi(w) \overline {\psi(-i \ol z)} d \mu(w) d \mu(z) \\
=&  \int_{V \times V} \lambda \exp (\frac{1}{4} (\ol z^t,  \ol w^t) \mc H(g) \arr{\ol z \\  \ol w})
\phi(w) \ol \psi(i z) d \mu(w) d \mu(z) \\
=& \Lambda(\lambda, g) \int_{V \times V}  \exp (\frac{1}{4} (\ol z^t,  \ol w^t) \mc H(g) \arr{\ol z \\  \ol w})
\phi(w) \ol \psi(i z) d \mu(w) d \mu(z)
\end{split}
\end{equation}
For any function $f \in \mc P_{2n}$, we define a function on $\mc S_{2n}$
\begin{equation}
W(f)(s)=\int_{V \times V}
\exp (\frac{1}{4} (\ol z^t, \ol w^t) s \arr{ \ol z \\ \ol w}) f(z, w) d \mu(z) d \mu(w).
\end{equation}
If we identify $\mc P_n \otimes \mc P_n$ with $\mc P_{2n}$ by
$$j: \phi(w) \otimes \psi(z) \rightarrow \phi(w) \ol \psi(iz),$$
then $W(j(\phi \otimes \psi))$ is the push-forward
 of $M_{\omega}(\phi \otimes \psi)(\lambda,g)/\Lambda(\lambda, g)$.

\begin{lem} 
Let $U(2n)$ act on $\mb C^{2n}$ as unitary operators. Let $U(2n)$ act on $\mc S_{2n}$
canonically. Then
$ {W}$ is a $U(2n)$ equivariant linear transform from
$\mc P_{2n}$ to $\mc O_{\mc S_{2n}}$. 
\end{lem}
Proof: For $g \in U(2n)$, $g$ acts on $\mc S_{2n}$ by $\tau(g) s=g s g^t$.
We compute
\begin{equation}
\begin{split}
{W} f (g^{-1} s (g^{-1})^t) &=\int_{\mb C^{2n}} \exp
(\frac{1}{4}(\ol z^t, \ol w^t) g^{-1} s (g^{-1})^t \arr{\ol z \\ \ol w}) f(\arr{z \\ w}) d \mu(z) d \mu(w) \\
&=\int_{\mb C^{2n}} \exp
(\frac{1}{4}(\ol z^t, \ol w^t) s \arr{\ol z \\ \ol w}) f(g^{-1} \arr{z \\ w}) d \mu(z) d \mu(w)
\end{split}
\end{equation}
Thus ${W}$ is $U(2n)$-equivariant. 
Q.E.D. 
 \\
\\
Let $\mc M_{1}$ be the linear span of matrix coefficients of 
$(\omega, \mc P_n)$. So $\mc M_1 \subset C^{\infty}(\meta)$. Let $\mc W_1$ to be the linear span of
$ W(f)$ with $f \in \mc P_{2n}$. Here $1$ is used to indicate that only one oscillator representation is involved as we will be discussing tensor products of the oscillator representation later. 
Notice that $\mc W_{1}$ is just the push-forward of $\mc M_{1}/ \Lambda$.
Let $\mc P^m_{2n}$ be the polynomials of $2n$ variables of homogeneous degree $m$.
Then it is well-known that $\mc P^m_{2n}$ is an irreducible unitary representation of $U(2n)$
with highest weight $(0,0,\ldots,0,-m)$. We obtain
\begin{lem}
$\mc W_1$ is contained in $\oplus_{i \in \mb N} C_{(0,0,\ldots,0,-i)}(\mc S_{2n})$.
\end{lem}
Observe that for an odd function f in $\mc P_{2n}$, 
$W f(s)=0$.
One can easily see that for any even $m$,
$${W} (\mc P_{2n}^m) \neq 0$$
since $(\omega(1) z^{\alpha}, w^{\alpha}) \neq 0$.
Therefore
\begin{thm}~\label{main}
$$\mc W_1=\oplus_{i \in \mb N} C_{(0,0,\ldots,0,-2i)}(\mc S_{2n})$$
In particular, $\mc W_1$ is a subspace of algebraic functions on $\mc S_{2n}$. Furthermore, $\mc W_1$ is the push-forward of $\mc M_1/\Lambda$ under $\mc H$.
\end{thm}
Here we regard $\mc M_1/\Lambda$ as a set of functions on $\xin$.
\subsection{The Contragradient $\omega^*$}
Let $\omega^*$ be the contragradient representation of $\omega$. So for $\delta \in \mc F_n^*$ and $u \in \mc F_n$, $\omega^*(g) \delta (u)= \delta(\omega(g^{-1} u))$. By Riesz representation theorem, $\delta$ can be identified with $u_{\delta}  \in \mc F_n$ such that
$$\delta(u)=(u,u_{\delta}).$$
It follows that $(u, u_{\omega^*(g) \delta}) = \omega^*(g) \delta (u)=\delta(\omega(g^{-1})u)=(\omega(g^{-1}) u, u_{\delta})=(u, \omega(g) u_{\delta})$. Therefore
$$ u_{\omega^*(g) \delta}= \omega(g) u_{\delta}.$$
Let $\mc P_n^c$ be $\mc P_n$ with complex conjugate linear scalar multiplication.
In other words, multiplication of $\lambda \in \mb C$ on $f \in \mc P_n^c$ is just $\ol \lambda f$. Then we may define a new representation $\omega^c$ by equipping
$(\omega, \mc F_n)$ with the complex conjugate linear scalar multiplication and the inner product
$$(u, v)_{\omega^c}=(v, u)_{\omega}.$$ Clearly, $\omega(g)= \omega^c(g)$. 
The matrix coefficient
$$(\omega^c(g) u, v)=(v, \omega(g) u)= \overline{(\omega(g)u, v)}= \overline{M_{\omega}(u \otimes v)(g)}.$$
Observe that
$\omega^* \cong \omega^c$. The identification is given by $\delta \rightarrow u_{\delta}$.
\begin{thm}
The push-forward of $M_{\omega^c}(u \otimes v)/\overline{\Lambda}$ spans $\overline{\mc W_1}=\oplus_{i \in \mathbb N} C_{(2i,0, \ldots, 0)} (\mc S_{2n})$.
\end{thm}

\subsection{Estimates on Matrix Coefficients}
Notice that functions in $\mc M_1$ are all bounded functions. Theorem ~\ref{main} implies 
\begin{lem}
For $u, v \in \mathcal P_n$, $| M_{\omega}(u \otimes v) | \leq C_{u,v} | \Lambda(\lambda, g)|$ for some constant $C_{u,v}$.
\end{lem}
Here $|\Lambda(\lambda, g)|=\prod_{i=1}^n (\sech H_i(g))^{\frac{1}{2}}$. For simplicity, write it as $| \Lambda(g)|$.
To obtain more precise estimates on $M_{\omega}$ for $f \in \mathcal P_n \otimes \mathcal P_n$ and beyond, we prove the following lemma.
\begin{lem}
Suppose that  $f \in L^1(z,w, \exp -\frac{1}{4}(\|z \|^2+ \|w \|^2) d z d w)$. Then
$$\| W f (s)\|_{max} \leq \| f \|_{L^1(z,w, \exp -\frac{1}{4}(\|z \|^2+ \|w \|^2) d z d w)}.$$
\end{lem}
Proof: Notice that 
$$\exp(\frac{1}{4}(\ol z^t, \ol w^t) s \arr{\ol z \\ \ol w}) \leq \exp(\frac{1}{4} (\|z \|^2+ \|w\|^2)).$$
We have $\forall \ s \in \mc S_{2n}$
\begin{equation}
\begin{split}
 | Wf (s)| = & |\int_{\mb C^{2n}} \exp
(\frac{1}{4}(\ol z^t, \ol w^t) s \arr{\ol z \\ \ol w}) f( \arr{z \\ w}) d \mu(z) d \mu(w)| \\
\leq & \int_{\mb C^{2n}} \exp
(\frac{1}{4}(\|z \|^2+ \|w\|^2)) |f( \arr{z \\ w})| \exp(-\frac{1}{2} (\|z \|^2+ \|w\|^2)) 
d z d w \\
= & \int_{\mb C^{2n}} |f( \arr{z \\ w})|  \exp -\frac{1}{4}(\|z \|^2+ \|w \|^2) d z d w) \\
= & \| f \|_{L^1(z,w, \exp -\frac{1}{4}(\|z \|^2+ \|w \|^2) d z d w)} 
\end{split}
\end{equation}
Q.E.D.

\begin{thm}~\label{es} Suppose that $u$ and $v$ are smooth vectors in $(\omega, \mc F_n)$. Then there exists a constant $C_n$ such that
$$ \| W (u \otimes v) (s)\|_{max} \leq C_n^2 \| u(z) v(w) (1+\|z\|)^{n+1} (1+ \|w\|)^{n+1} \|_{L^2(z,w, \exp -\frac{1}{2} (\|z\|^2 + \| w\|^2) d z d w)}$$ and for every $g \in \meta$
$$\| M_{\omega}(u \otimes v)(g) \| \leq  C_n^2 \| u(z) v(w) (1+\|z\|)^{n+1} (1+ \|w\|)^{n+1} \|_{L^2(z,w, \exp -\frac{1}{2} (\|z\|^2 + \| w\|^2) d z d w)} | \Lambda(g) | $$
\end{thm}
Proof: Since $u$ is smooth, $ z^{\alpha} u \in \mc F_n$. By Schwartz inequality, 
$$( \int | u(z)| \exp -\frac{1}{4}\|z \|^2  d z)^2 \leq \int |u(z) (1+ \|z\|)^{n+1} \exp -\frac{1}{4} \|z \|^2 |^2  d z
\int |(1+ \|z\|)^{-n-1}|^2 d z$$
Therefore, there exists a constant $C_n$, only depending on $n$, such that
$$ \|u\|_{L^1(z, \exp -\frac{1}{4}(\|z\|^2) d z)} \leq C_n \| u(z) (1+\|z\|)^{n+1} \|_{L^2(z, \exp -\frac{1}{2} \|z\|^2 d z)}.$$
By the previous lemma,
$$ | W (u \otimes v)(s)| \leq C_n^2 \| u(z) v(w) (1+\|z\|)^{n+1} (1+ \|w\|)^{n+1} \|_{L^2(z,w, \exp -\frac{1}{2} (\|z\|^2 + \| w\|^2) d z d w )}.$$
It follows that
$$ | M_{\omega}( u \otimes v)(g) | \leq  C_n^2 \| u(z) v(w) (1+\|z\|)^{n+1} (1+ \|w\|)^{n+1} \|_{L^2(z,w, \exp -\frac{1}{2} (\|z\|^2 + \| w\|^2) d z d w )} |\Lambda(g)|  .$$

\subsection{Tensor Products of Oscillator Representations}
Let 
$$\omega_{p,q}=\overbrace{\omega \otimes \omega \otimes \ldots \otimes \omega}^p \otimes
\overbrace{\omega^c \otimes \omega^c \otimes \ldots \otimes \omega^c}^q,$$
$$\mc P_{p,q}= (\otimes^{p} \mc P_n) \otimes (\otimes^q \mc P_n^c).$$
Let $\mc M_{p,q}$ be the linear span of  matrix coefficients of $(\omega_{p,q}, \mc P_{p,q})$. Then $\mc M_{p,q}/ (\Lambda^p \overline{\Lambda}^q) \subset C^{\infty}(\xin)$. 
Define 
$$W_{p,q}= \otimes^p W \otimes \otimes^q \overline{W}:
f_1 \otimes f_2 \otimes \cdots \otimes f_{p+q} \in \otimes^{p} \mc P_{2n} \otimes \otimes^q \mc P_{2n}^c
  \rightarrow \prod_{i=1}^p
W(f_i) \prod_{i=p+1}^{p+q} \overline{W(f_i)}.$$
Let
$\mc W_{p,q}$ be the image of $W_{p,q}$.
According to a Theorem of Kashiwara-Vergne (~\cite{kv}), the $U(2n)$-types in 
$(\otimes^{p} \mc P_{2n}) \otimes (\otimes^q \mc P_{2n}^c)$
are parametrized by 
$$(m_1 \geq m_2 \geq \ldots \geq m_s \geq \overbrace{0,0,\ldots 0}^{2n-s-t} \geq n_1 \geq n_2 \geq \ldots \geq n_t)$$
with $s \leq q$, $t \leq p$ and $s+t \leq 2n$. It is easy to see that the set of parameters
is exactly $S_{q,p}$.
According to Helgason's theorem,
the $U(2n)$ types in $\mc O_{\mc S_{2n}}$ are parametrized by even highest weights.
We have
\begin{lem}~\label{mpq}
$\mc W_{p,q} = \oplus_{\lambda \in S_{q,p} \cap 2 \mathbb Z^{2n}} C_{\lambda}(\mc S_{2n})$.
\end{lem}
Proof: We only need to show that every $C_{\lambda}(\mc S_{2n})$ with $\lambda \in S_{q,p} \cap 2 \mathbb Z^{2n}$ indeed occurs in
$\mc W_{p,q}$. Let 
$$\lambda= (m_1 \geq m_2 \geq \ldots \geq m_s \geq \overbrace{0,0,\ldots 0}^{2n-s-t} \geq n_1 \geq n_2 \geq \ldots \geq n_t)$$
with $ s \leq q, t \leq p$.
Recall that 
$$\mc W_1=\oplus_{i \in \mb N} C_{(0,0,\ldots,0,-2i)}(\mc S_{2n})$$
Let $v_{-2i}$ be the highest weight vector of $C_{(0,0, \ldots,0,-2i)}(\mc S_{2n})$.
Then $v_{\lambda}=\prod_{i=1}^t v_{n_i} \prod_{i=1}^s \overline{v_{-m_i}}$ is a highest weight vector
in $C_{\lambda}(\mc S_{2n})$. Furthermore,
$v_{\lambda}$ is in $\mc W_{p,q}$ since 
$$\mc W_{p,q}=span\{\prod_{i=1}^{p} f_i \prod_{i=p+1}^{p+q} \overline{f_i},
 \mid f_i \in \mc W_1 \ \ \forall \, \, i \}.$$
Thus $C_{\lambda}(\mc S_{2n}) \subset \mc W_{p,q}$. Q.E.D. \\
\\
By Theorem ~\ref{main}, we have
\begin{thm}~\label{push}
The matrix coefficients of $\omega_{p,q}$ can be written as
$\Lambda^p \overline{\Lambda}^q f(\mc H(g))$
with $f \in \mc W_{p,q}$. So $\mc W_{p,q}$ is the push forward of $\mc M_{p,q}/ \Lambda^p \overline{\Lambda^q}$ under $\mc H$.
\end{thm}

\subsection{The function ${\det}_{12}$}
Recall
$${\det}_{12}(s)=(\wedge^n(s)  e_1 \wedge e_2 \wedge \ldots \wedge e_n, e_{n+1} \wedge e_{n+1}
\wedge \ldots \wedge e_{2n}).$$
Let $V_{\lambda}$ be an irreducible constituent of the linear span of
$$ \{ \det(g s g^t) \in C^{\infty}(\mc S_{2n}) \mid  g \in U(2n) \}.$$  Here $\lambda$ is a dominant weight.
Then
\begin{itemize}
\item $\lambda$ is even.
\item $V_{\lambda}$ contains a vector of weight $(-1,-1,-1, \ldots, -1)$.
\item $\lambda$ is in the convex hull spanned by permutations of $(0, 0, \ldots, 0, \overbrace{-2,-2, \ldots,-2 }^n)$.
\end{itemize}
Thus $\lambda=(0,0, \ldots, 0, \overbrace{-2,-2, \ldots,-2 }^n)$.
\begin{lem}~\label{phi}
The function ${\det}_{12}(s) \in C_{\lambda}(\mc S_{2n})$ with $\lambda={(0,0, \ldots, 0, \overbrace{-2,-2, \ldots,-2 }^n)}$. Furthermore, for every $s \in \mc S_{2n}$,
$${\overline{{\det}_{12}(s)}}= (-1)^n {{\det}}^{-1} (s) {{\det}_{12}(s)}.$$
\end{lem}
Proof: We only need to prove the second statement. Consider the generalized Cartan decomposition
$$s= \arr{k_1 & 0 \\ 0 & k_2} \arr{ \cos \theta & - i \sin \theta \\ -i \sin \theta & \cos \theta} \arr{k_1^t & 0 \\ 0 & k_2^t}.$$
The determinant
$$ \overline{{\det}_{12}(s)}= \det(i \overline{k_1} \sin \theta \overline{k_2})=
(-1)^n (\det k_1)^{-2} (\det k_2)^{-2} \det (-i k_1 \sin \theta k_2) =(-1)^n \det s^{-1} {\det}_{12}(s).$$
Q.E.D.\\
\\
The function ${\det}^{-1}$ is regarded as a function on $\mc S_{2n}$. It is of the weight $(2,2 \ldots,2 )$ due to the left multiplication by $g^{-1}$.

\subsection{An Orthogonal decomposition of $L^2_{-}(\meta)$}
\begin{thm}~\label{decmeta}
$(\mc M_{p,q}, \mc M_{p_1,q_1})=0$ if $p_1+q_1=p+q=2n+1$, $p_1 \neq p$ and
$p_1 \equiv p \pmod 2$. Furthermore
$$L^2_{-}(\meta= \hat{\oplus}_{p+q=2n+1, p \ odd} \mc M_{p,q}.$$
$$L^2_{-}(\meta= \hat{\oplus}_{p+q=2n+1, p \ even} \mc M_{p,q}.$$
 \end{thm}
Proof:  Suppose that $p_1+q_1=p+q=2n+1$, $p_1 \neq p$ and
$p_1 \equiv p \pmod 2$. Assume that $p$ is odd and $q $ is even. For $f \in \mc M_{p,q}$, recall that
$I(f)=(f \Lambda^{-2n-1})^0.$ According to Theorem ~\ref{iso}, it suffices to show that
$$(I(\mc M_{p,q}), I(\mc M_{p_1,q_1}))=0$$
and
$$L^2(\mc S_{2n})= \hat{\oplus}_{p+q=2n+1, p \ odd} I(\mc M_{p,q}).$$
According to
Theorem ~\ref{phi12} and Lemma ~\ref{phi}, the push-forward of $(\frac{\overline{\Lambda
}}{\Lambda})^q$
is just 
$$(-1)^{\frac{nq}{2}}(\frac{\overline{{\det}_{12}(s)}}{{\det}_{12}(s)})^{\frac{q}{2}}=
(-1)^{n +\frac{nq}{2}} {\det} (s)^{\frac{-q}{2}} .$$ 
 It follows that 
$$I(f)=(f \Lambda^{-2n-1})^0=(f \Lambda^{-p} \overline{\Lambda}^{-q})^0( ( \frac{\overline{\Lambda
}}{\Lambda})^q)^0 .$$
By Theorem ~\ref{push},  $(f \Lambda^{-p} \overline{\Lambda}^{-q})^0 \in \mc W_{p,q}$. So
$$I(\mc M_{p,q})= {\det} (s)^{-\frac{q}{2}} \mc W_{p,q}.$$
From Theorem ~\ref{mpq},
$$I(\mc M_{p,q})=  \oplus_{\lambda \in S_{q,p}} {\det} s^{\frac{-q}{2}} C_{\lambda}(\mc S_{2n})= \oplus_{\lambda \in \mathbf{q}+ S_{q,p}} C_{\lambda}(\mc S_{2n}).$$
Our first assertion follows from Theorem ~\ref{deco}. Our second assertion follows similarly from Theorem ~\ref{deco}. Q.E.D.\\
\\
We denote the closure of $\mc M_{p,q}$ by $cl(\mc M_{p,q})$. Then
$$L^2_{-}(\meta= \oplus_{p+q=2n+1, p \ odd} cl(\mc M_{p,q});$$
$$L^2_{-}(\meta= \oplus_{p+q=2n+1, p \ even} cl(\mc M_{p,q}).$$
\begin{cor}[Orthonormal Basis]
There exists, $B_{p,2n+1-p}$, a set of elements in 
$$\mc P_{p, 2n+1-p} \otimes \mc P_{p, 2n+1-p}  $$
 such that
$\{ M_{\omega_{p, 2n+1-p}}(b) \}_{b \in B_{p, 2n+1-p}}$ is an orthonormal basis in $cl(\mc M_{p, 2n+1-p})$.
In particular, all $M_{\omega_{p, 2n+1-p}}(b)$ is bounded by a multiple of $\| \Lambda(g) \|^{2n+1}$.
\end{cor}
This Corollary follows from Theorem ~\ref{push}.

\section{Signature of Genuine Tempered Representations}

\begin{thm} Let $p+q=2n+1$. Let $\meta$ act on $L^2_{-}(\meta)$ from both left and right. Then $cl(\mc M_{p,q})$ is a left and right subrepresentation.  Let $P_{p,q}$ be the projection of $L^2_{-}(\meta)$ onto $cl(\mc M_{p,q})$.
Then $P_{p,q}$ is $\meta$-equivariant.
\end{thm}
Proof: Fix $(p,q)$. Let $B$ be the linear space spanned by
$$\omega_{p,q} (g) u   \qquad g \in \meta, u \in  \mc P_{p,q}.$$ 
Clearly, $\meta$ acts on $B$ from left. Let $\mc M_B$ be the set of  matrix coefficients with respect to $B$. Then $\mc M_{B}$ is a $\meta$-space. 
{\bf Claim:} $ \mc M_{p,q} \subseteq \mc M_B \subseteq cl(\mc M_{p,q})$.\\
\\
Clearly $\mc M_{p,q} \subseteq \mc M_{B}$. To prove the second inclusion, let $u, v \in \mc P_{p,q}$. Fix $g_1 \in \meta$. It suffices to show that $g \rightarrow (\omega_{p,q}(g g_1) u, v)$ is in $cl(\mc M_{p,q})$. Since $\omega_{p,q}(g_1) u$ is an analytic function, $\omega_{p,q}(g_1) u$ has a Taylor expansion. Let $u_i$ be its $i$-th Taylor polynomial in terms of the total degree. Since $\omega_{p,q} (g_1) u$ is smooth,  $z^{\alpha} w^{\beta} u_i (z) v(w) \rightarrow z^{\alpha} w^{\beta} \omega_{p,q}(g_1) u (z) v (w)$ in $L^2(z,w, d \mu(z,w))$. By Theorem ~\ref{es},
$$ W_{p,q} (u_i \otimes v)(s)  \rightarrow W_{p,q} (\omega_{p,q}(g g_1) u \otimes v)(s)$$
in sup-norm in $C(\mc S_{2n})$.  In particular, $W_{p,q} (u_i \otimes v)(s)$ is uniformly bounded by a constant function. Recall that 
$$| M_{\omega_{p,q}}(u_i \otimes v)(g)|=|W_{p,q}(u_i \otimes v)(\mc H(g))| |\Lambda(g)\|^{2n+1}$$
and $|\Lambda(g)|^{2n+1} \in L^2(\xin)$.
It follows that $ M_{\omega_{p,q}} (u_i \otimes v)(g)  \rightarrow M_{\omega_{p,q}}(\omega_{p,q}(g g_1) u \otimes v)(g)$ in $L^2(\meta)$. Thus,
$g \rightarrow (\omega_{p,q}(g g_1) u, v)$ is in $cl(\mc M_{p,q})$.\\
\\
Now $g \in \meta$ preserves on $\mc M_B$ and $\mc M_B$ is dense in $cl(\mc M_{p,q})$. By the unitarity of $\omega_{p,q}(g)$,  $g$ preserves on $cl(\mc M_B) = cl(\mc M_{p,q})$. So
$cl(\mc M_{p,q})$ is a subrepresentation of $L^2_{-}(\meta)$ under the left regular action.
The rest of the theorem follows easily. Q.E.D.

\subsection{Signature of Genuine Discrete Series Representation}
Let $D$ be a genuine discrete series representation of  $\meta$. Let $u,v$ be $K$-finite vectors in $D$.
Then $D_{u,v}(g)=(D(g)u,v)$ is in $L^2_{-}(\meta)$. Now under the decomposition 
$$L^2_{-}(\meta)= {\oplus}_{p+q=2n+1, p \ odd} cl(\mc M_{p,q}),$$
$D_{u,v}(g)$ decomposes into a sum of functions in $cl(\mc M_{p,q})$. What we want to show is that
$D_{u,v}(g)$ in actually in one and only one of $cl(\mc M_{p,q})$.
\begin{thm}~\label{discretesig}
Fix a parity of $p$ and let $p+q=2n+1$. Let $D$ be a genuine discrete series representation of  $\meta$. Then there exists
a unique $p$ such that the matrix coefficients of $D$ are in $cl(\mc M_{p,q})$. Equivalently, there exists a unique $p$ such that $D$ is a subrepresentation of $cl(\mc M_{p,q})$.
\end{thm}
Proof: Let $D$ be a genuine discrete series representation. Let $H_D \hat{\otimes} H_D^*$ be the irreducible subrepresentation of $\meta \times \meta$ occurring discretely in $L^{2}_{-}(\meta)$.  Consider $P_{p,q}( H_D \hat{\otimes} H_D^*)$. There exists a $(p,q)$ such that $P_{p,q}( H_D \hat{\otimes} H_D^*) \neq 0$. Since $P_{p,q}$ is $\meta \times \meta$-equivariant, $P_{p,q}( H_D \hat{\otimes} H_D^*)$ is a $\meta \times \meta$ representation. Because of the irreducibility of $ H_D \hat{\otimes} H_D^*$, $P_{p,q}( H_D \hat{\otimes} H_D^*)$ is equivalent to $H_D \hat{\otimes} H_D^*$. By Schur's Lemma, $P_{p,q}|_{ H_D \hat{\otimes} H_D^*}= \lambda I_{ H_D \hat{\otimes} H_D^*}$. Since
$P_{p,q}$ is a projection, $P_{p,q}|_{ H_D \hat{\otimes} H_D^*}=I_{ H_D \hat{\otimes} H_D^*}$. So for any other
$p_1+q_1=2n+1$ and $p_1 \equiv p \pmod 2$, $P_{p_1,q_1}( H_D \hat{\otimes} H_D^*)=0$. Therefore, there exists a unique
$p$ such that $H_D \hat{\otimes} H_D \subset cl(\mc M_{p,q})$. Our assertion follows. Q.E.D.
\begin{defn}~\label{discrete}
Let $D$ be a genuine discrete series representation of  $\meta$. Fix a parity of $p$. We call $(p, 2n+1-p)$ the signature of $D$ for $p$ odd (even) if the matrix coefficients of $D$ are in $cl(\mc M_{p,2n+1-p})$. We denote the odd $p$ by $p_o$ and the even $p$ by $p_e$. We call
$(p_e, 2n+1-p_e; p_o, 2n+1-p_o)$ the signature of $D$.
\end{defn}
\subsection{Signature of Genuine Tempered Representations}
Let $p+q=2n+1$.
Recall that the matrix coefficients of $\omega_{p,q}$ are in $L^2(\meta)$. By a Theorem of Cowling-Haagerup-Howe ~\cite{chh},
the representations in the weak closure of $\omega_{p,q}$ are all tempered. Let
$\Pi_{p,q}^2(\meta)$ be the support of $\omega_{p,q}$. Let $\Pi^2_{-}(\meta)$ be the set of equivalence classes of irreducible tempered representations with $\pi(\epsilon)=-1$. The set $\Pi_{-}^2(\meta)$ is a discrete union of
$\mathbb R^{n}/W$ where $n$ is a nonnegative integer and $W$ is a certain Weyl group. The Fell topology on $\Pi_{-}^2(\meta)$ is the natural topology. The representation $L^2_{-}(\meta)$ is supported on $\Pi^2_{-}(\meta)$ (see 14.12 ~\cite{wallach}).
\begin{thm}~\label{temsig} Fix a parity of $p$.
$\Pi^2_{-}(\meta)$ is the disjoint union of $\Pi_{p,2n+1-p}^2(\meta)$.
\end{thm}
Proof: Fix a parity of $p$. First of all, we want to prove that $\Pi^2_{-}(\meta)=\cup_{p} \Pi_{p,2n+1-p}^2(\meta)$. 
Put $U= \Pi^2_{-}(\meta)- \cup_{p} \Pi_{p,2n+1-p}^2(\meta)$. Then $U$ is open with respect to the Fell topology.
If the Plancherel measure of $U$ is zero, then $L^2_{-}(\meta)$ is supported on the closure of $\cup_{p} \Pi_{p,2n+1-p}^2(\meta)$. Since $\Pi_{p,2n+1-p}^2(\meta)$ is already closed, we have $$\Pi^2_{-}(\meta)= \cup_{p} \Pi_{p,2n+1-p}^2(\meta).$$
We are done. Now suppose that the Plancherel measure of $U$ is nonzero. Let $L^2_{-}(\meta, U_1)$ be a nontrivial subspace of $L^2_{-}(\meta)$ supported on a closed subset $U_1$ in $U$. Let $f$ be a nonzero function in $L^2_{-}(\meta, U_1)$. Notice that $\omega_{p,2n+1-p}$ is supported on $\Pi_{p,2n+1-p}^2(\meta)$. It follows that $f \perp \mc M_{p, 2n+1-p}$ for all $p$ with a fixed parity. 
 By Theorem ~\ref{decmeta}, $f(g)=0$ for all $g$. This is a contradiction. So $U=\emptyset$ and
 $$\Pi^2_{-}(\meta)=\cup_{p} \Pi_{p,2n+1-p}^2(\meta).$$\\
 \\
 Next, we want to show that if $p_1 \neq p_2$ and $p_1 \equiv p_2 \pmod 2$ then $\Pi_{p_1, 2n+1-p_1}^2(\meta) \cap \Pi_{p_2, 2n+1-p_2}^2(\meta) = \emptyset$. Notice that Theorem ~\ref{decmeta} only implies that the intersection is of Plancherel measure zero. Suppose that $(\pi, H_{\pi}) \in \Pi_{p_1, 2n+1-p_1}^2(\meta) \cap \Pi_{p_2, 2n+1-p_2}^2(\meta)$. Let $u$ be a nonzero $K$-finite vector in $H_{\pi}$. Suppose that
 $u$ is in the $K$-type $\sigma$. Consider $M_{\pi}(u \otimes u)$. Since $\pi$ is weakly contained in
 $\Pi_{p_1, 2n+1-p_1}^2(\meta)$, there exists a sequence of elements $u_i \in \mc P_{p_1, 2n+1-p_1}$ such that $M_{\omega_{p_1,2n+1-p_1}}(u_i \otimes u_i) \rightarrow M_{\pi}(u \times u)$ uniformly over any compact set. We may further assume that $u_i$ are all of the $K$-type $\sigma$ and $\{ \| u_i\| \}$ is bounded.
 Let $p \neq p_1$ and $ p \equiv p_1 \pmod 2$.
 By Theorem ~\ref{decmeta}, $M_{\omega_{p_1,2n+1-p_1}}(u_i \otimes u_i) \perp \mc M_{p, 2n+1-p}$. Now we want to apply the Dominated Convergence Theorem to show that $M_{\pi}(u \otimes u) \perp \mc M_{p, 2n+1-p}$.
 By a Theorem of Cowling-Haagerup-Howe, $| M_{\omega_{p_1,2n+1-p_1}}(u_i \otimes u_i)| $ are uniformly bounded by $C \Xi(g)$ where $\Xi(g)$ is Harish-Chandra's basic spherical function ~\cite{chh}. All functions in $\mc M_{p, 2n+1-p}$ are bounded by $|\Lambda|^{2n+1}$. It is easy to see that $\Xi(g) |\Lambda|^{2n+1} \in L^1(\meta)$. So by the Dominated Convergence Theorem, $M_{\pi}(u \otimes u) \perp \mc M_{p, 2n+1-p}$ for all $p \neq p_1$ and $p \equiv p_1 \pmod 2$. Similarly, $M_{\pi}(u \otimes u) \perp \mc M_{p, 2n+1-p}$ for all $p \neq p_2$ and $p \equiv p_1 \pmod 2$. Since $p_1 \neq p_2$, we have  $M_{\pi}(u \otimes u) \perp \mc M_{p, 2n+1-p}$ for all $p \equiv p_1 \pmod 2$. By Theorem ~\ref{decmeta}, $M_{\pi}(u \otimes u)=0$. We reach a contradiction. We have finished showing that
 $$\Pi_{p_1, 2n+1-p_1}^2(\meta) \cap \Pi_{p_2, 2n+1-p_2}^2(\meta) = \emptyset.$$
 
\begin{defn}~\label{tempered}
Let $\pi$ be a genuine tempered representation of  $\meta$. Fix a parity of $p$. We call $(p, 2n+1-p)$ the signature of $\pi$ for $p$ odd (even) if $\pi \in \Pi_{p,2n+1-p}^2(\meta)$. We denote the odd $p$ by $p_o$ and the even $p$ by $p_e$. We call
$(p_e, 2n+1-p_e; p_o, 2n+1-p_o)$ the signature of $\pi$.
\end{defn}
From the proof of Theorem ~\ref{temsig}, we obtain
\begin{cor}~\label{24} Fix a parity of $p$. $\pi \in \Pi_{-}^2(\meta)$ is in $\Pi_{p, 2n+1-p}^2(\meta)$ if and only if
the $K$-finite matrix coefficients of $\pi$ are all perpendicular to $\mc M_{p_1, 2n+1-p_1}$ for every $p_1 \neq p$ with $p_1 \equiv p \pmod 2$.
\end{cor}
Notice that any $K$-finite matrix coefficient $f(g)$ is bounded by a multiple of $\Xi(g)$ and functions in $\mc M_{p_1,2n+1-p_1}$ are all bounded by $\Lambda^{2n+1}$. So it makes sense to say that $$(f(g), \mc M_{p_1,2n+1-p_1})=0.$$ In fact, we can make a stronger statement.
\begin{cor}~\label{pq} Fix a parity of $p$.  $\pi \in \Pi_{-}^2(\meta)$ is in $\Pi_{p, 2n+1-p}^2(\meta)$ if and only if there exists a
 $K$-finite matrix coefficient $f$ of $\pi$ and a matrix coefficient $F$ in $\mc M_{p, 2n+1-p}$ such that
 $$\int_{\meta} f(g) F(g) d g \neq 0.$$
\end{cor}
Proof: We prove the only if part by contradiction. Suppose $(f(g), \mc M_{p, 2n+1-p})=0$ for every $K$-finite matrix coefficient of $f$. Since $\pi \in \Pi_{p, 2n+1-p}^2(\meta)$, by the previous Corollary, $(f(g), \mc M_{p_1,2n+1-p_1})=0$ for every $p_1 \neq p$ with $p_1 \equiv p \pmod 2$.
By Theorem ~\ref{decmeta}, $f =0$ almost everywhere. $\pi$ is not a representation. \\
\\
Now suppose that there exists a
 $K$-finite matrix coefficient $f$ of $\pi$ and a matrix coefficient $F$ in $\mc M_{p, 2n+1-p}$ such that
 $$\int_{\meta} f(g) F(g) d g \neq 0.$$
 If $\pi \notin \Pi_{p, 2n+1-p}^2(\meta)$, then by Theorem ~\ref{temsig}, there exists a $p_1 \neq p$ and $p_1 \equiv p \pmod 2$ such that $\pi \in \Pi_{p_1, 2n+1-p_1}^2(\meta)$. By Cor. \!~\ref{24}, $(f, \mc M_{p, 2n+1-p})=0$. So $\int_{\meta} f(g) F(g) d g =0.$ This is a contradiction. Q.E.D.

\begin{thm} Let $D$ be a genuine discrete series representation of $\meta$.
Fix a parity of $p$. Then there is a unique $p$ such that $D$ is equivalent to a subrepresentation of $\omega_{p,2n+1-p}$.
\end{thm}
Proof:  Since $D$ is tempered, $D$ is in the weak closure of $\omega_{p, 2n+1-p}$ for a unique $p$. Because
$D$ is isolated in $\Pi_{-}^2(\meta)$, $D$ must occur as a subrepresentation of $\omega_{p, 2n+1-p}$. Q.E.D.\\
\\
Comparing this theorem with Theorem ~\ref{discretesig}, we see that Definition ~\ref{tempered} coincides with
Definition ~\ref{discrete} for discrete series representations. 
\subsection{ Howe Duality and A Theorem of Adam and Barbasch}
To define signature beyond tempered representations, we encounter some serious technical difficulties. The trouble is that the push forward of the matrix coefficients of a nontempered $\pi$ may fail to be locally integrable near the boundary with respect to a certain measure. At this stage, we do not know how to overcome this obstacle. Nevertheless, we can bypass this by relating our results to the Howe duality (~\cite{howe}). \\
\\
Let $\Pi(G)$ denote the admissible dual of a semisimple Lie group $G$. Let $\mc R(\meta, \omega_{p, 2n+1-p})$ be those $\pi$ in $\Pi(\meta)$ such that $\pi$ occurs as a quotient of $\omega_{p, 2n+1-p}$ infinitesimally. Howe's theorem then states that there is a one-to-one correspondence between $\mc R(\meta, \omega_{p, 2n+1-p})$ and $\mc R(O(p, 2n+1-p), \omega_{p, 2n+1-p})$.
We denote the Howe duality for $(\meta, O(p, 2n+1-p))$ by $\theta(2n; p, 2n+1-p)$ and $\theta(p, 2n+1-p;2n)$. \\
\\
A result of Adam and Barbasch states that $\theta(2n; p, 2n+1-p)$ induces a one-to-one correspondence between $\Pi_{-}(\meta)$ and the disjoint union of $\Pi(SO(p, 2n+1-p))$ for $p$ with a fixed parity. In particular, $\Pi_{-}(\meta)$ is the disjoint union of $\mc R(\meta, \omega_{p, 2n+1-p})$ for $p$ with a fixed parity. We can now define a signature of an irreducible genuine representation $\pi$ to be $(p, 2n+1-p)$ if $\pi \in \mc R(\meta, \omega_{p, 2n+1-p})$.
Is this definition different from ours when $\pi$ is tempered? The answer is No.\\
\\
We shall now use the theory developed in ~\cite{theta} to show that
\begin{thm}~\label{closure} Let $\pi$ be a tempered representation. $\pi \in \Pi_{p, 2n+1-p}^2(\meta)$ if and only if $\pi \in \mathcal R(\meta, \omega_{p,2n+1-p})$ (see ~\cite{howe}).
\end{thm}
Proof: Let $\pi$ be a genuine tempered representation of $\meta$. Suppose that $\pi \in \Pi_{p, 2n+1-p}^2(\meta)$. Then by Definition 3.2.1 ~\cite{basic}, $\pi$ is in the semistable range of
$\theta(2n;p, 2n+1-p)$. Let $f$ be a $K$-finite matrix coefficient of $\pi$. By Cor ~\ref{pq}, $(f(g), \mc M_{p, 2n+1-p}) \neq 0$. By Theorem 1.1 ~\cite{theta}, $\pi \in \mathcal R(\meta, \omega_{p,2n+1-p})$. \\
\\
Conversely, since $\pi$ is in the semistable range of
$\theta(2n;p, 2n+1-p)$, by Theorem 1.1 ~\cite{theta}, $(f(g), \mc M_{p, 2n+1-p}) \neq 0$ for a $K$-finite matrix coefficient $f(g)$ of $\pi$. By Cor ~\ref{pq},
$\pi \in \Pi_{p, 2n+1-p}^2(\meta)$. 
Q.E.D.\\
\\
Finally, let me state a theorem concerning tempered representations of $SO(p,q)$ and $\meta$.
\begin{thm} Fix a parity of $p$.  $\theta$ induces a one-to-one correspondence between
$\Pi^2_{-}(\meta)$ and $\cup \Pi^2(SO(p, 2n+1-p))$.
\end{thm}
Proof: Fix a parity of $p$. We give our proof in two steps.\\
\\
First, let $\pi$ be a tempered irreducible representation of $SO(p, 2n+1-p)$. By Def. 3.2.1 (~\cite{basic}), $\pi$ is in the semistable range of $\theta(p, 2n+1-p; 2n)$. By Theorem ~\cite{non}, $\theta(p,2n+1-p;2n)(\pi)$ exists. By Lemma 6.2.1 and Example 2 (~\cite{basic}), $\theta(p,2n+1-p;2n)(\pi)$ is a tempered irreducible representation.\\
\\
Second, let $\pi \in \Pi_{-}^2 (\meta)$. Let $(p, 2n+1-p)$ be its signature, depending on the parity of $p$. By Def. 3.2.1 (~\cite{basic}), $\pi$ is in the semistable range of $\theta(2n; p, 2n+1-p)$. From theorem ~\ref{closure}, $\theta(2n; p, 2n+1-p)(\pi)$ is nonvanishing. By Lemma 6.3.1 and Example 3 (~\cite{basic}), $\theta(2n; p, 2n+1-p)(\pi)$ is a tempered irreducible representation. Q.E.D. \\
\\
We formulate the following conjecture concerning Howe duality (see ~\cite{howe}).
\begin{conj}
Let $(G_1, G_2)$ be a real reductive dual pair. Then $\pi$ occurs
in $\mc R(MG_1, \omega)$ if and only if the matrix coefficients of $\pi$ can be approximated
by the restrictions of matrix coefficients of $\omega$ onto $MG_1$ uniformly on compacta.
\end{conj}
\subsection{$p+q=2n+2$: Some Results of Moeglin}
C. Moeglin treats the case $p+q=2n+2$ with $p,q$ both even. One of her results is that every $\pi \in \Pi(\xin)$ occurs as a quotient of
$\omega_{p,2n+2-p}$ for some even integer $p$. In this section, I shall briefly show that  every tempered representation $\pi$ occurs as a quotient of $\omega_{p,2n+2-p}$ for an even integer $p^e$ and an {\bf odd integer} $p^o$. We need the following lemma.
\begin{lem}
 The linear span of $\{ \mc M_{p, 2n+2-p} \mid p \ even \}$ is dense in $L^2(\xin)$ and the linear span of $\{ \mc M_{p, 2n+2-p} \mid p \ odd \}$ is also dense in $L^2(\xin)$.
 \end{lem}
This Lemma follows directly from Theorem ~\ref{push} and Lemma ~\ref{mpq}. The proof is omitted. We now have
 \begin{thm}
 Fix a parity of $p$. Let $\pi \in \Pi^2(\xin)$. Then $\pi \in \mc R(\xin, \omega(p, 2n+2-p; 2n))$ for at least one $p$.
 \end{thm}
Finally, let me make a conjecture concerned with the question about the growth of matrix coefficients at $\infty$. 
\begin{conj}
Let $\pi$ be a unitary irreducible representation of $\xin$ with an
integrable infinitesimal character. Let $M_{\pi}(u \otimes v)(g)$ be a $K$-finite matrix coefficient.
Then the push forward  $M_{\pi}(u \otimes v)(g)^0 $ onto $\mc S_{2n}$ is analytic. 
\end{conj}

\end{document}